# FITTING AN ERROR DISTRIBUTION IN SOME HETEROSCEDASTIC TIME SERIES MODELS[1]

By Hira L. Koul and Shiqing Ling

*Michigan State University and Hong Kong University of Science and Technology*

This paper addresses the problem of fitting a known distribution to the innovation distribution in a class of stationary and ergodic time series models. The asymptotic null distribution of the usual Kolmogorov–Smirnov test based on the residuals generally depends on the underlying model parameters and the error distribution. To overcome the dependence on the underlying model parameters, we propose that tests be based on a vector of certain weighted residual empirical processes. Under the null hypothesis and under minimal moment conditions, this vector of processes is shown to converge weakly to a vector of independent copies of a Gaussian process whose covariance function depends only on the fitted distribution and not on the model. Under certain local alternatives, the proposed test is shown to have nontrivial asymptotic power. The Monte Carlo critical values of this test are tabulated when fitting standard normal and double exponential distributions. The results obtained are shown to be applicable to GARCH and ARMA–GARCH models, the often used models in econometrics and finance. A simulation study shows that the test has satisfactory size and power for finite samples at these models. The paper also contains an asymptotic uniform expansion result for a general weighted residual empirical process useful in heteroscedastic models under minimal moment conditions, a result of independent interest.

**1. Introduction.** Let $\{y_i : i \in \mathbb{Z} := 0, \pm 1, \pm 2, \ldots\}$ be a strictly stationary and ergodic real time series. Often the finite-dimensional distributions of such series are characterized by the stationary distribution and the conditional distribution of $y_i$, given the past. One problem of interest is to fit this

Received November 2003; revised February 2005.
[1]Supported by the Hong Kong RGC Grants HKUST4765/03H and HKUST6022/05P.
*AMS 2000 subject classifications.* Primary 62F05, 62M10; secondary 60G10.
*Key words and phrases.* Nonlinear time series models, goodness-of-fit test, weighted empirical process.

conditional distribution. In general, this is a difficult problem. However, in some special time series models where this conditional distribution is determined by the innovation distribution, it is possible to obtain reasonable answers. In particular, in this paper we shall focus on the generalized autoregressive conditionally heteroscedastic (GARCH) and ARMA–GARCH models.

To describe these models, let $q, r$ be known positive integers and $\Theta_1(\Theta_2)$ be a subset of the $q(r)$-dimensional Euclidean space $\mathbb{R}^q(\mathbb{R}^r)$, and let $\Theta = \Theta_1 \times \Theta_2$. In the models of interest one observes the process $y_i$ such that, for some sequences of past measurable functions $\mu_i$ from $\Theta_1$ to $\mathbb{R}$ and $h_i$ from $\Theta$ to $\mathbb{R}^+ := (0, \infty)$, and for some $\theta' = (\theta_1', \theta_2')$, $\theta_1 \in \Theta_1$, $\theta_2 \in \Theta_2$,

$$(1.1) \qquad \eta_i := \frac{y_i - \mu_i(\theta_1)}{\sqrt{h_i(\theta)}}, \qquad i \in \mathbb{Z},$$

are independently and identically distributed (i.i.d.) standardized r.v.'s. Here, "past measurable" means that, for every $s := (s_1, s_2) \in \Theta, s_1 \in \Theta_1$, the functions $\mu_i(s_1)$ and $h_i(s)$ should be $\mathcal{F}_{i-1}$ measurable, where $\mathcal{F}_i$ is the $\sigma$-field generated by $\{\eta_i, \eta_{i-1}, \ldots, y_0, y_{-1}, \ldots\}$, $i \in \mathbb{Z}$. Let $F$ denote the common distribution function (d.f.) of the errors $\{\eta_i\}$, and $F_0$ be a known d.f. The problem of fitting the conditional distribution of $y_i$, given $\mathcal{F}_{i-1}$, in the model (1.1) is now equivalent to testing the goodness-of-fit hypothesis

$$H_0 : F = F_0 \quad \text{vs.} \quad H_1 : F \neq F_0.$$

The knowledge of the error distribution is important in statistics, in particular, in value at risk (VaR). In economics and finance, VaR is a single number measuring the risk of a financial position. For example, when $y_i$ is a process of daily returns, the VaR for a one-day horizon of a portfolio is the 95th conditional quantile of the distribution of $y_{i+1}$, given the information available at time $i$. After estimating the parameter $\theta$, the VaR for a one-day position of $y_i$ and probability 0.05 is $\mu_i(\theta_1) - 1.6449\sqrt{h_i(\theta)}$, provided the $\eta_i$'s have the standard normal distribution and the estimated parameter $\theta$ is correct. This means that, with probability 0.95, the potential loss of holding that position the next day is $\mu_i(\theta_1) - 1.6449 \times \sqrt{h_i(\theta)}$. Clearly, the knowledge of the error distribution plays a crucial role in determining this probability, and, hence, in evaluating VaR via model (1.1). Thus, it is important to test the hypothesis $H_0$ in practice. For more on VaR, see, for example, [29].

The goodness-of-fit testing problem under the i.i.d. setup has a long history; see, for example, a collection of papers in [7], and references therein. A commonly used test is based on the Kolmogorov–Smirnov statistic. The primary reason for this is that this test is distribution free, that is, the null distribution of this statistic does not depend on $F_0$. However, when the i.i.d.

sequence such as $\eta_i$ in the model (1.1) is not observed and has to be estimated from a special model, this property of being distribution free is lost even asymptotically. These kinds of problems have been extensively investigated in the literature in various models, including regression models; see, for example, Durbin [9], Loynes [28] and Koul [14], among others.

In the context of time series models, Boldin [2] observed that for the zero mean linear autoregressive (AR) models the tests based on the residual empirical process are asymptotically distribution free (ADF) for fitting an error d.f. as long as it has zero mean, finite variance and bounded second derivative. The condition of bounded second derivative was later relaxed in [13] to requiring having only a uniformly continuous density, while the condition of having zero mean is crucial for the validity of the ADF property. Boldin [3] and Koul [14] observed that similar facts hold for moving-average (MA) models. However, this is no longer true of many nonlinear time series models such as threshold AR models as noted in [15], nonstationary AR models in [21], ARCH/GARCH models in [4, 5, 11], [16], Chapter 8 and [1]. For ARCH/GARCH models, Horváth, Kokoszka and Teyssière [12] proposed a parametric boostrap method for testing Gaussianity of the errors, but only studied its validity via Monte Carlo experiments. Some of the results of Kulperger and Yu [18] can be used to test for the error moments but not for the error d.f.

The main difficulty is that the asymptotic null distribution of the empirical process of the residuals in the heteroscedastic time series model (1.1) depends not only on the distribution $F_0$, but also on the model functions $\mu_i$ and $h_i$. Under certain regularity conditions, the dependence on $\mu_i$ can be eliminated by using a vector of certain weighted residual empirical processes, as is shown here. Roughly speaking, the weights in these processes are asymptotically orthogonal to the space generated by the slopes of the functions $\mu_i$. Moreover, under $H_0$ this vector of processes is shown to converge weakly to a vector of independent copies of a Gaussian process. The covariance function of this Gaussian process depends only on $F_0$ and is the same as that of the one arising when fitting a d.f. up to an unknown scale parameter. The Monte Carlo critical values when $F_0$ is either the standard normal or double exponential d.f. are tabulated. Under certain local alternatives, the proposed test is shown to have nontrivial asymptotic power.

Section 2 describes the test statistic and the main results. Section 3 shows that our assumptions are naturally satisfied by a class of GARCH and ARMA–GARCH models. Some simulation results and an application of the proposed test to the Hang Seng Index in the Hong Kong stock market are also given in Section 3. Section 4 contains some proofs. Theorem 1.1 of [17] has been found very useful in many nonlinear homoscedastic time series models and has been extended for nonstationary time series in [27]. Theorem 4.1 below gives its extension to heteroscedastic time series models and is used to prove the main results of Section 2.

**2. Test statistics and the main results.** This section describes the proposed test statistics and their asymptotic behavior.

Let $p$ be a positive integer and $Y_0$ denote either the vector $(y_0, y_{-1}, \ldots, y_{1-p})'$ or the vector $(y_0, y_{-1}, \ldots)'$. In either case $Y_0$ is chosen to be independent of $\eta_i, i \geq 1$. Let $\{y_1, \ldots, y_n\}$ be observations obeying the model (1.1). We shall assume the following:

(2.1) The functions $\mu_i, h_i$ are twice continuously differentiable, for all $i$.

Let $\varepsilon_i(s_1) = y_i - \mu_i(s_1)$, $\eta_i(s) := \varepsilon_i(s_1)/\sqrt{h_i(s)}$, $1 \leq i \leq n$. For any differentiable function $g$ on $\Theta$, let $\dot{g}$ denote its differential. Thus, for example, $\dot{\mu}_i(s_1) = \partial \mu_i(s_1)/\partial s_1$, $s_1 \in \Theta_1$. Moreover, for an $s \in \Theta$, we write

$$\dot{\mu}_i(s) = \begin{pmatrix} \dot{\mu}_i(s_1) \\ 0 \end{pmatrix}, \qquad \dot{h}_i(s) = \begin{pmatrix} \dot{h}_{1i}(s) \\ \dot{h}_{2i}(s) \end{pmatrix},$$

$$\dot{h}_{ki}(s) := \frac{\partial h_i(s)}{\partial s_k}, \qquad k = 1, 2, \ 1 \leq i \leq n.$$

Thus, $\dot{\mu}_i(s)$ and $\dot{\mu}_i(s_1)$ are $(q+r) \times 1$ and $q \times 1$ vectors, respectively. Similarly, $\dot{h}_{1i}(s)$ and $\dot{h}_{2i}(s)$ are $q \times 1$ and $r \times 1$ vectors, respectively. We denote the true parameter by $\theta$ and $\eta_i = \eta_i(\theta)$.

We need the following additional assumptions:

(2.2) $F_0$ has an absolutely continuous density $f_0$ with $\sup_{x \in \mathbb{R}} |x| f_0(x) < \infty$ and has positive and finite Fisher information for location and scale, that is,

(2.3)
$$0 < \int (\dot{f}_0/f_0)^2 \, dF_0 < \infty, \qquad 0 < \int \left(1 + \frac{x \dot{f}_0(x)}{f_0(x)}\right)^2 dF_0(x) < \infty.$$

$$E\psi_0(\eta_1)\varphi_0(\eta_1) = 0,$$

where $\psi_0 := \dot{f}_0/f_0$ and $\varphi_0(x) := (1 + x\psi_0(x))/2$, $x \in \mathbb{R}$.

For example, under (2.2), (2.3) holds whenever $F_0$ is symmetric around zero. Let

$$b_1 := E\psi_0^2(\eta_1), \qquad b_2 := E\varphi_0^2(\eta_1), \qquad \mathcal{B} := E\begin{pmatrix} b_1 & 0 \\ 0 & b_2 \end{pmatrix},$$

$$W_i(s) := \begin{pmatrix} \dot{\mu}_i(s_1)/h_i^{1/2}(s) & \dot{h}_{1i}(s)/h_i(s) \\ 0 & \dot{h}_{2i}(s)/h_i(s) \end{pmatrix} = \begin{pmatrix} W_{11,i}(s) & W_{12,i}(s) \\ 0 & W_{22,i}(s) \end{pmatrix}, \quad \text{say,}$$

$$\mathcal{I}(\theta) := E_\theta W_1(\theta) \mathcal{B} W_1(\theta)'.$$

Under $H_0$, $\mathcal{B}$ is known and $\mathcal{I}(\theta)$ is the Fisher information matrix appropriate in the present setup. We shall assume that, under $H_0$, $\mathcal{I}(\theta)$ is positive definite. Let $\hat{\theta}_n$ be an estimator of $\theta$ satisfying

$$(2.4) \qquad n^{1/2}(\hat{\theta}_n - \theta) = -\mathcal{I}(\theta)^{-1} n^{-1/2} \sum_{i=1}^{n} W_i(\theta) \begin{pmatrix} \psi_0(\eta_i) \\ \varphi_0(\eta_i) \end{pmatrix} + o_p(1).$$

We are now ready to introduce the needed weighted empirical processes,

$$K_n(x,s) := \frac{1}{2\sqrt{n}} \sum_{i=1}^{n} W_{22,i}(s)[I(\eta_i(s) \leq x) - F_0(x)], \qquad x \in \mathbb{R}, s \in \Theta.$$

Note that this is an $r \times 1$ vector of processes. The tests of $H_0$ will be based on $K_n(x, \hat{\theta}_n)$. For example, the Kolmogorov–Smirnov test statistic is

$$K_n = \sup_{x \in \mathbb{R}} [K_n(x, \hat{\theta}_n)' \widehat{\mathcal{I}}_n^{-1} K_n(x, \hat{\theta}_n)], \qquad \widehat{\mathcal{I}}_n := \frac{1}{4n} \sum_{i=1}^{n} W_{22,i}(\hat{\theta}_n) W_{22,i}(\hat{\theta}_n)'.$$

The following theorem is useful in obtaining the limiting distribution of the tests based on the process $K_n(\cdot, \hat{\theta}_n)$. In it the matrix of the second derivative of $h_i$ is denoted by $\ddot{h}_i$ and the matrix norm is the Euclidean norm.

THEOREM 2.1. *Suppose the model* (1.1) *and* (2.1)–(2.4) *hold. In addition, suppose the following holds. For an open neighborhood $U_\theta$ of $\theta$,*

$$(2.5) \qquad E_\theta \sup_{s \in U_\theta} \frac{\|\dot{\mu}_1(s)\|^2}{h_1(\theta)} < \infty, \qquad E_\theta \sup_{s \in U_\theta} \frac{\|\dot{h}_1(s)\|^2}{h_1(s)h_1(\theta)} < \infty.$$

$$(2.6) \qquad E_\theta \sup_{s \in U_\theta} \frac{\|\ddot{h}_1(s)\|}{h_1(s)} < \infty, \qquad E_\theta \sup_{s \in U_\theta} \frac{\|\dot{h}_1(s)\|^2}{h_1^2(s)} < \infty.$$

*Then, for every $0 < b < \infty$,*

$$K_n(x, \theta + n^{-1/2}t)$$

$$(2.7) \qquad = K_n(x, \theta) + \frac{1}{2} n^{-1} \sum_{i=1}^{n} W_{22,i}(\theta) \left[ \frac{\dot{\mu}_i(\theta)'}{\sqrt{h_i(\theta)}} f_0(x) + x f_0(x) \frac{1}{2} \frac{\dot{h}_i(\theta)'}{h_i(\theta)} \right] t$$

$$+ u_p(1),$$

*where $u_p(1)$ is a sequence of stochastic processes tending to 0 uniformly in $x \in \mathbb{R}$ and $\|t\| \leq b$, in probability.*

A proof of this theorem is sketched in Section 4. Here we shall now illustrate its usefulness in obtaining the limiting null distribution of $K_n$.

Now, consider the coefficient of $t/2$ in the second term on the right-hand side of (2.7). By the ergodic theorem, uniformly in $x \in \mathbb{R}$, it converges to

$$(2.8) \qquad E_\theta W_{22,1}(\theta) \left[ \frac{\dot{\mu}_1(\theta)'}{\sqrt{h_1(\theta)}} f_0(x) + x f_0(x) \frac{1}{2} \frac{\dot{h}_1(\theta)'}{h_1(\theta)} \right].$$

To reduce the effect of the location related parameters, namely, of $\theta_1$ and $\mu_1(\theta_1)$, on this expression, it suffices to assume

$$E_\theta\{[W_{11,1}(\theta), W_{12,1}(\theta)]W_{22,1}(\theta)'\} = 0. \tag{2.9}$$

Under this assumption, the expression in (2.8) is equivalent to

$$\tfrac{1}{2}xf_0(x)[0', E_\theta W_{22,1}(\theta)W_{22,1}(\theta)'] = \tfrac{1}{2}xf_0(x)[0', H'_\theta], \qquad \text{say.}$$

Next, let $I_{11}(\theta) := b_1 E_\theta W_{11,1}(\theta)W_{11,1}(\theta)' + b_2 E_\theta W_{21,1}(\theta)W_{21,1}(\theta)'$. Under (2.9),

$$\mathcal{I}(\theta) = \begin{pmatrix} W_{21,1}(\theta)' & 0 \\ I_{11}(\theta) & 0 \\ 0 & b_2 H_\theta \end{pmatrix}.$$

Moreover, the first term on the right-hand side of the expansion (2.4) now becomes

$$-n^{-1/2} \sum_{i=1}^n \begin{pmatrix} I_{11}^{-1}(\theta)(W_{11,i}\psi_0(\eta_i) + W_{12,i}(\theta)\varphi_0(\eta_i)) \\ (b_2 H_\theta)^{-1} W_{22,i}(\theta)\varphi_0(\eta_i) \end{pmatrix}.$$

Upon combining these facts with (2.7), we obtain that, uniformly in $x \in \mathbb{R}$,

$$K_n(x, \hat\theta_n) = K_n(x, \theta) - \frac{xf_0(x)}{4b_2\sqrt{n}} \sum_{i=1}^n W_{22,i}(\theta)\varphi_0(\eta_i) + o_p(1)$$

$$= \frac{1}{2\sqrt{n}} \sum_{i=1}^n W_{22,i}(\theta)\left[I(\eta_i \le x) - F_0(x) - \frac{xf_0(x)}{2b_2}\varphi_0(\eta_i)\right] + o_p(1).$$

Let $Z_n(x)$ denote the vector of the leading process on the right-hand side above, and

$$\rho(x,y) := \left[F_0(x \wedge y) - F_0(x)F_0(y) - \frac{1}{4b_2}xy f_0(x)f_0(y)\right], \qquad x, y \in \mathbb{R}.$$

Let $Z$ be a vector of $r$ independent mean zero Gaussian processes, with $\text{Cov}(Z(x), Z(y)) := \rho(x,y)I_{r\times r}$. Using a conditioning argument and the above weak convergence result, one readily obtains that $\text{Cov}(Z_n(x), Z_n(y)) = \frac{1}{4}\rho(x,y)H_\theta$, and that $2H_\theta^{-1/2}K_n(x, \hat\theta_n) \Longrightarrow Z(x)$. This fact, together with the fact that $\|\hat{\mathcal{I}}_n - H_\theta\| = o_p(1)$, yields the following corollary.

COROLLARY 2.1. *Suppose the assumptions of Theorem* 2.1 *and* $H_0$ *hold. If in addition* (2.9) *holds, then* $K_n \xrightarrow{d} \sup_{x \in \mathbb{R}} \|Z(x)\| =: K$.

Observe that this limiting distribution depends only on the error d.f. $F_0$. Thus, its critical values can be approximated by simulation. When $\eta_i \sim N(0,1)$, the corresponding critical values are given in Table 1. This table is constructed as follows. First, we approximate $\{Z(x) : x \in [-4.0, 4.0]\}$ by $\{Z(x_i) : i = 1, \ldots, 2000\}$ with $x_{i+1} - x_i = 8/2000$. The distribution of $K$ is then approximated by $\sup_{1 \leq i \leq 2000} \|Z(x_i)\|$. We use 10,000 independent replications to obtain the percentages of $K$. By increasing the ranges of $x$ and $i$ and the number of replications, the percentages of $K$ have only ignorable differences from those in Table 1. In the case $f_0(x) = e^{-|x|}/2, x \in \mathbb{R}$, the corresponding critical values are in Table 2. The simulation method is the same as that for Table 1 except now the range of $x$ is taken to be $[-8.0, 8.0]$, since this density has a longer tail than the standard normal density. For other distributions, the critical values can be obtained through a similar method.

We now indicate the behavior of the asymptotic power of this test under the local alternatives $H_{1n} : F_n = (1 - n^{-1/2}\delta)F_0 + n^{-1/2}\delta\tilde{F}$, where $0 < \delta < 1$ and $\tilde{F}$ is a d.f. Assume that $\{\eta_i; i \geq 1\}$ are independent of $Y_0$ under $H_{1n}$. Let $P_{0n}$ and $P_{1n}$ be the joint distributions of $(y_1, \ldots, y_n)$ under $H_0$ and $H_{1n}$, respectively. In Section 5 we show that $P_{0n}$ and $P_{1n}$ are contiguous in the sense of Le Cam [19]. The following theorem gives the asymptotic power of the test statistic $K_n$ against these local alternatives.

THEOREM 2.2. *In addition to the conditions of Theorem* 2.1, *assume that $\tilde{F}$ has zero mean and finite variance. Then, under $H_{1n}$, $K_n(x, \hat{\theta}_n) \Longrightarrow$*

TABLE 1
*Upper percentage points of $K$ with $F_0 = N(0, 1)$*

| $\alpha \backslash r$ | 1 | 2 | 3 | 4 | 5 | 6 | 7 | 8 | 9 | 10 |
|---|---|---|---|---|---|---|---|---|---|---|
| 0.01 | 2.465 | 3.150 | 3.737 | 4.361 | 4.769 | 5.173 | 5.642 | 6.065 | 6.508 | 6.938 |
| 0.03 | 1.890 | 2.595 | 3.100 | 3.640 | 4.094 | 4.468 | 4.946 | 5.336 | 5.677 | 6.083 |
| 0.05 | 1.650 | 2.289 | 2.804 | 3.267 | 3.751 | 4.118 | 4.553 | 4.922 | 5.298 | 5.679 |
| 0.10 | 1.317 | 1.891 | 2.382 | 2.822 | 3.262 | 3.628 | 4.033 | 4.384 | 4.716 | 5.106 |
| 0.15 | 1.113 | 1.666 | 2.126 | 2.552 | 2.980 | 3.325 | 3.685 | 4.022 | 4.367 | 4.720 |
| 0.20 | 0.988 | 1.498 | 1.931 | 2.339 | 2.750 | 3.089 | 3.436 | 3.768 | 4.101 | 4.445 |

TABLE 2
*Upper percentage points of $K$ with $F_0 =$ double exponential distribution*

| $\alpha \backslash r$ | 1 | 2 | 3 | 4 | 5 | 6 | 7 | 8 | 9 | 10 |
|---|---|---|---|---|---|---|---|---|---|---|
| 0.01 | 2.402 | 3.149 | 3.702 | 4.298 | 4.809 | 5.173 | 5.683 | 5.937 | 6.360 | 6.845 |
| 0.03 | 1.876 | 2.523 | 3.073 | 3.569 | 4.015 | 4.399 | 4.872 | 5.260 | 5.611 | 6.015 |
| 0.05 | 1.630 | 2.250 | 2.781 | 3.218 | 3.680 | 4.067 | 4.500 | 4.865 | 5.247 | 5.607 |
| 0.10 | 1.299 | 1.873 | 2.344 | 2.788 | 3.222 | 3.605 | 3.969 | 4.335 | 4.680 | 5.051 |
| 0.15 | 1.098 | 1.640 | 2.092 | 2.533 | 2.933 | 3.279 | 3.639 | 3.991 | 4.317 | 4.691 |
| 0.20 | 0.961 | 1.464 | 1.902 | 2.314 | 2.703 | 3.046 | 3.399 | 3.729 | 4.065 | 4.046 |

$2H_\theta^{1/2}\tilde{Z}(x)$, where $\tilde{Z}(x)$ is an $r$-vector of Gaussian processes with mean vector $m_\theta(x) := E_\theta W_{22,1}(\theta)(F(x) - \tilde{F}_0(x))$, and covariance matrix equal to that of the process $Z(x)$ in Corollary 2.1.

This theorem shows that the test statistic $K_n$ has nontrivial asymptotic power under the specified local alternatives, provided $m_\theta(x) \neq 0$, for some $x \in \mathbb{R}$.

**3. Some applications.** This section verifies the applicability of the results in the previous section to the GARCH and ARMA–GARCH models. These models are among the most important models in the field of time series, econometrics and finance.

3.1. *GARCH models.* This subsection considers Bollorslev's GARCH $(p_1, p_2)$ model defined by the equations

$$(3.1) \qquad y_i = \eta_i \sqrt{h_i} \quad \text{and} \quad h_i = \alpha_0 + \sum_{j=1}^{p_1} \alpha_j y_{j-i}^2 + \sum_{j=1}^{p_2} \beta_j h_{j-i},$$

where $\eta_i$ are standardized i.i.d. r.v.'s. Clearly, it is an example of the model (1.1) with $\mu_i \equiv 0 = q$, $h_i$ as given above, and $s = s_2 = (\alpha_0, \alpha_1, \ldots, \alpha_{p_1}, \beta_1, \ldots, \beta_{p_2})'$ is the unknown parameter vector and its true value is $\theta = \theta_2 = (\alpha_{00}, \alpha_{01}, \ldots, \alpha_{0p_1}, \beta_{01}, \ldots, \beta_{0p_2})'$. The parameter space is $\Theta = \{s : \sum_{j=1}^{p_1} \alpha_j + \sum_{j=1}^{p_2} \beta_j \leq \rho_0, \underline{\beta}_0 \leq \alpha_0 \leq \tilde{\beta}_0, \text{ and } \underline{\beta} \leq \alpha_i, \beta_j \leq \tilde{\beta}, i = 1, \ldots, p_1, j = 1, \ldots, p_2\}$ for some constants $\rho_0 \in (0,1)$, $0 < \underline{\beta}_0 \leq \tilde{\beta}_0$ and $0 < \underline{\beta} \leq \tilde{\beta} < 1$. Note that $\Theta$ is compact. Assume that $\theta$ is an interior point of $\Theta$ and $\sum_{j=1}^{p_1} \alpha_j z^j$ and $1 - \sum_{j=1}^{p_2} \beta_j z^j$ have no common zeros for each $s \in \Theta$. Bougerol and Picard [6] have shown that, under these conditions, the model (3.1) is strictly stationary and ergodic with $E_\theta y_0^2 < \infty$. We now verify the remaining assumptions of the previous section. Let

$$\left(1 - \sum_{j=1}^{p_2} \beta_j B^j\right)^{-1} \left(\sum_{j=1}^{p_1} \alpha_j B^j\right) = \sum_{j=1}^{\infty} \beta_{0j}(s) B^j,$$

$$\left(1 - \sum_{j=1}^{p_2} \beta_j B^j\right)^{-1} = \sum_{j=1}^{\infty} \beta_{1j}(s) B^j,$$

where $B$ denotes the backward shift operator. Then from (3.1) one has

$$h_i(s) = \left(1 - \sum_{j=1}^{p_2} \beta_j B^j\right)^{-1} \left(\sum_{j=1}^{p_2} \alpha_j B^j\right) y_i^2 = \sum_{j=1}^{\infty} \beta_{0j}(s) y_{i-j}^2, \qquad \text{say,}$$

$$\frac{\partial h_i(s)}{\partial \alpha_0} = 1 + \sum_{j=1}^{p_2} \beta_j \frac{\partial h_{i-j}(s)}{\partial \alpha_0} = \left(1 - \sum_{j=1}^{p_2} \beta_j\right)^{-1},$$

$$\frac{\partial h_i(s)}{\partial \alpha_k} = y_{i-k}^2 + \sum_{j=1}^{p_2} \beta_j \frac{\partial h_{i-j}(s)}{\partial \alpha_k} = \left(1 - \sum_{j=1}^{p_2} \beta_j B^j\right)^{-1} y_{i-k}^2,$$

$$\frac{\partial h_i(s)}{\partial \beta_k} = h_{i-k}(s) + \sum_{j=1}^{p_2} \beta_j \frac{\partial h_{i-j}(s)}{\partial \beta_k} = \left(1 - \sum_{j=1}^{p_2} \beta_j B^j\right)^{-1} h_{i-k}(s),$$

$$\frac{\partial^2 h_i(s)}{\partial \alpha_0^2} = 0,$$

$$\frac{\partial^2 h_i(s)}{\partial \beta_k^2} = 2\frac{\partial h_i(s)}{\partial \beta_k} + \sum_{j=1}^{p_2} \beta_j \frac{\partial^2 h_{i-j}(s)}{\partial \beta_k^2} = 2\left(1 - \sum_{j=1}^{p_2} \beta_j B^j\right)^{-1} \frac{\partial h_i(s)}{\partial \beta_k}.$$

Similar expressions can be obtained for $\partial^2 h_i(s)/\partial \alpha_k^2$ and $\partial^2 h_i(s)/\partial \alpha_k \partial \beta_{k'}$, $1 \leq k \leq p_1, 1 \leq k' \leq p_2$. From these expressions, it follows that the above $\{h_i(s)\}$ satisfy (2.1). From Lemma 2.1 in [22] and the compactness of $\Theta$, it follows that, for some $0 < \rho < 1$, $\sup_{s \in \Theta} \beta_{0j}(s) = O(\rho^j)$ as $j \to \infty$.

Next, note that here $W_{11,i}(s) \equiv 0, W_{12,i}(s) \equiv 0$ and

$$W_{22,i}(s) = \frac{1}{h_i(s)} \left[\frac{\partial h_i(s)}{\partial \alpha_0}, \frac{\partial h_i(s)}{\partial \alpha_1}, \ldots, \frac{\partial h_i(s)}{\partial \alpha_{p_1}}, \frac{\partial h_i(s)}{\partial \beta_1}, \ldots, \frac{\partial h_i(s)}{\partial \beta_{p_2}}\right]'$$

$$= \frac{1}{h_i(s)} \frac{\partial h_i(s)}{\partial s}.$$

Thus, (2.9) is a priori satisfied. Further, assume that $E\eta_1^4 < \infty$. Francq and Zakoïan [10] show that the quasi-MLE is $\sqrt{n}$-consistent. Using this quasi-MLE as an initial estimator and one-step iteration as in (3.1) in [23], we can obtain a new estimator $\hat{\theta}_n$. Ling [23] showed that this $\hat{\theta}_n$ satisfies the condition (2.4) under (2.2).

Next, we verify (2.5)–(2.6). Because $\mu_i \equiv 0$, the first part of (2.5) is a priori satisfied. By Lemma A.2 in [24], there are a neighborhood $U_\theta$ of $\theta$ and a finite universal constant $C$ and a $\rho \in (0,1)$ such that, for any $\gamma > 0$, and for all $i$,

$$\max\left(\sup_{s \in U_\theta} \left\|\frac{1}{h_i(s)} \frac{\partial^2 h_i(s)}{\partial s \, \partial s'}\right\|, \sup_{s \in U_\theta} \left\|\frac{1}{h_i(s)} \frac{\partial h_i(s)}{\partial s}\right\|\right) \leq C\left(1 + \sum_{j=1}^{\infty} \rho^j |y_{i-j}|\right)^\gamma.$$

Thus, (2.6) holds. By Lemma A.1(ii) and A.5 in [24], there exists a constant $\tilde{\iota} \in (0,1)$ such that $\sup_{s \in U_\theta}[h_i(s)/h_i(\theta)] \leq C(1 + \sum_{j=1}^{\infty} \rho^j |\varepsilon_{i-j}|)^{2-2\tilde{\iota}}$. Taking

$\gamma \leq \tilde{\iota}$, it follows that

$$E_\theta \sup_{s \in U_\theta} \frac{\|\dot{h}_1(s)\|^2}{h_1(s)h_1(\theta)} < \frac{C}{\alpha_0} E_\theta \sup_{s \in U_\theta} \left[ \frac{h_1(s)}{h_1(\theta)} \left( 1 + \sum_{j=1}^\infty \rho^j |y_{-j}| \right)^{2\gamma} \right] < \infty.$$

Thus, (2.5) holds. This completes the verification of all the assumptions of Corollary 2.1.

3.2. ARMA(1, 1)–GARCH(1, 1) *models*. Consider the ARMA(1, 1)–GARCH(1, 1) model defined by the equations

$$
\begin{aligned}
y_i &= a y_{i-1} + b \varepsilon_{i-1} + \varepsilon_i, \\
\varepsilon_i &= \eta_i \sqrt{h_i} \quad \text{and} \quad h_i = \alpha_0 + \alpha \varepsilon_{i-1}^2 + \beta h_{i-1},
\end{aligned}
\tag{3.2}
$$

where $\eta_i$ is a sequence of i.i.d. standardized r.v.'s and symmetric around zero. We shall assume the usual condition: $\alpha, \beta \in (0, 1)$, $a, b \in (-1, 1)$, and $a + b \neq 0$. This model is an example of the model (1.1) with $q = 2$, $r = 3$, $s_1 = (a, b)'$, $s_2 = (\alpha_0, \alpha, \beta)'$, $s = (s_1, s_2)'$, $\varepsilon_i(s_1) = y_i - a y_{i-1} - b \varepsilon_{i-1}(s_1)$, $\mu_i(s_1) \equiv y_i - \varepsilon_i(s_1)$, and $h_i(s) = \alpha_0 + \alpha \varepsilon_{i-1}^2(s_1) + \beta h_{i-1}(s)$. We shall assume that, for some positive numbers $\underline{a}_0, \tilde{a}_0, \underline{a}_1, \tilde{a}_1$, $\underline{a}_0 \leq \alpha_0 \leq \tilde{a}_0$, $\underline{a}_1 \leq \alpha \leq \tilde{a}_1$, $\underline{a}_1 \leq \beta \leq \tilde{a}_1$ and $\alpha + \beta \leq \rho_0$ with $\rho_0 \in (0, 1)$. The parameter space $\Theta$ is now the subset of $\mathbb{R}^5$ whose members satisfy these conditions and is a priori compact. The true parameter values $\theta' = (\theta_1', \theta_2')$ are taken to be in the interior of $\Theta$, where $\theta_1 = (a_0, b_0)'$, $\theta_2 = (\alpha_{00}, \alpha_0, \beta_0)'$. From [25], it follows that, under these conditions, the model (3.2) is strictly stationary, ergodic, and that $E_\theta y_0^2 < \infty$.

To verify the remaining assumptions of the previous section for this model, proceed as follows. From (3.2), one obtains

$$\frac{\varepsilon_i(s_1)}{\partial a} = -y_{i-1} + b \frac{\varepsilon_{i-1}(s_1)}{\partial a} = -\sum_{j=1}^\infty b^{j-1} y_{i-j},$$

$$\frac{\varepsilon_i(s_1)}{\partial b} = -\varepsilon_{i-1}(s_1) + b \frac{\varepsilon_{i-1}(s_1)}{\partial b} = -\sum_{j=1}^\infty b^{j-1} \varepsilon_{i-j}(s_1),$$

$$\frac{\partial h_i(s)}{\partial a} = 2\alpha \varepsilon_{i-1}(s_1) \frac{\varepsilon_{i-1}(s_1)}{\partial a} + \beta \frac{\partial h_{i-1}(s)}{\partial a} = 2\alpha \sum_{j=1}^\infty \beta^{j-1} \varepsilon_{i-j}(s_1) \frac{\varepsilon_{i-j}(s_1)}{\partial a},$$

$$\frac{\partial h_i(s)}{\partial b} = 2\alpha \varepsilon_{i-1}(s_1) \frac{\varepsilon_{i-1}(s_1)}{\partial b} + \beta \frac{\partial h_{i-1}(s)}{\partial b} = 2\alpha \sum_{j=1}^\infty \beta^{j-1} \varepsilon_{i-j}(s_1) \frac{\varepsilon_{i-j}(s_1)}{\partial b},$$

$$\frac{\partial h_i(s)}{\partial \alpha_0} = (1 - \beta)^{-1},$$

$$\frac{\partial h_i(s)}{\partial \alpha} = \varepsilon_{i-1}^2(s_1) + \beta \frac{\partial h_{i-1}(s)}{\partial \alpha} = \sum_{j=1}^{\infty} \beta^{j-1} \varepsilon_{i-j}^2(s_1),$$

$$\frac{\partial h_i(s)}{\partial \beta} = h_{i-1}(s) + \beta \frac{\partial h_{i-1}(s)}{\partial \beta} = \sum_{j=1}^{\infty} \beta^{j-1} h_{i-j}(s).$$

Thus, here

$$W_{11,i}(s) = \frac{1}{\sqrt{h_i(s)}} \left[ \frac{\partial \varepsilon_i(s_1)}{\partial a}, \frac{\partial \varepsilon_i(s_1)}{\partial b} \right]' = \frac{1}{\sqrt{h_i(s)}} \frac{\partial \mu_i(s)}{\partial s_1},$$

$$W_{21,i}(s) = \frac{1}{h_i(s)} \left[ \frac{\partial h_i(s)}{\partial a}, \frac{\partial h_i(s)}{\partial b} \right]' = \frac{1}{h_i(s)} \frac{\partial h_i(s)}{\partial s_1},$$

$$W_{22,i}(s) = \frac{1}{h_i(s)} \left[ \frac{\partial h_i(s)}{\partial \alpha_0}, \frac{\partial h_i(s)}{\partial \alpha}, \frac{\partial h_i(s)}{\partial \beta} \right]' = \frac{1}{h_i(s)} \frac{\partial h_i(s)}{\partial s_2}.$$

The symmetry of $F_0$ around zero implies (2.9) here. Further, assume that $E\eta_1^4 < \infty$. Ling [24] showed that the self-weighted quasi-MLE is $\sqrt{n}$-consistent. Using this estimator as an initial estimator and one-step iteration as in (3.1) in [23], we can obtain a new estimator $\hat{\theta}_n$. Ling [23] showed that this $\hat{\theta}_n$ satisfies condition (2.4) under (2.2) for more general fractional FARIMA–GARCH models, which includes (3.2) as a special case.

Next, we verify (2.5) and (2.6). Since $E\varepsilon_i^2(\theta_1) < \infty$, there exist some $0 < C < \infty$ and $0 < \rho < 1$ such that

$$(3.3) \quad E_\theta \sup_{s \in U_\theta} \frac{\|\dot{\mu}_1(s)\|^2}{h_1(\theta)} \leq E_\theta \sup_{s \in U_\theta} \|\dot{\mu}_1(s)\|^2 \leq CE_\theta \left( \sum_{k=1}^{\infty} \rho^k |\varepsilon_{i-k}| \right)^2 < \infty,$$

$$(3.4) \quad E_\theta \sup_{s \in U_\theta} \left\| \frac{1}{\sqrt{h_i(s)}} \frac{\partial h_i(s)}{\partial s_1} \right\|^2 \leq CE_\theta \left[ \sup_{s \in U_\theta} \sum_{j=1}^{\infty} \beta^{(j-1)/2} \left\| \frac{\varepsilon_{i-j}(s_1)}{\partial s_1} \right\| \right]^2$$

$$\leq CE_\theta \left[ \sum_{k=1}^{\infty} \rho^k |\varepsilon_{i-k}| \right]^2 < \infty.$$

By Lemma A.2 in [24], there are a neighborhood $U_\theta$ of $\theta$ and universal constants $0 < C < \infty$, $0 < \rho < 1$, such that, for any $\gamma > 0$, and uniformly in $i$ and $s$, $\|h_i^{-1}(s) \partial h_i(s)/\partial s_2\| \leq C(1 + \sum_{j=1}^{\infty} \rho^j |y_{i-j}|)^\gamma$. By Lemma A.1(ii) and A.5 in [24], there exists a constant $\tilde{\iota} \in (0,1)$ such that $\sup_{s \in U_\theta} [h_i(s)/h_i(\theta)] \leq C(1 + \sum_{j=1}^{\infty} \rho^j |\varepsilon_{i-j}|)^{2-2\tilde{\iota}}$. Thus, taking $\gamma \leq \tilde{\iota}$, it follows that

$$(3.5) \qquad E \sup_{s \in U_\theta} \frac{1}{h_i(\theta) h_i(s)} \left\| \frac{\partial h_i(s)}{\partial s_2} \right\|^2 < \infty.$$

From (3.3)–(3.5), we thus obtain that $E_\theta \sup_{s \in U_\theta}[\|\dot{h}_1(s)\|^2/(h_1(s)h_1(\theta))] \leq \alpha_{00}^{-1} E_\theta \sup_{s \in U_\theta} \|\dot{h}_1(s)\|^2/h_1(s) < \infty$, thereby verifying (2.5) in the present case. The condition (2.6) is verified similarly, thereby completing the verification of the assumptions of Corollary 2.1 in this example.

3.3. *Empirical results.* This section first examines the performance of the test statistic $K_n$ in finite samples through Monte Carlo experiments. The following AR(1)–GARCH(1, 1) model is used:

$$(3.6) \quad y_i = ay_{i-1} + \varepsilon_i, \qquad \varepsilon_i = \eta_i \sqrt{h_i}, \qquad h_i = \alpha_0 + \alpha \varepsilon_{i-1}^2 + \beta h_{i-1},$$

where the true parameters are $(a, \alpha_0, \alpha, \beta) = (0.5, 0.025, 0.25, 0.5)$ and $\eta_i$ are i.i.d. In the experiments, we take the sample sizes $n = 200$ and 400. One thousand replications are used. The null distribution of $\eta_i$ is standard normal and its alternatives are given as follows:

$A_1: \quad \sqrt{3/5}t_5, \qquad A_2: \quad \sqrt{1/2}t_4, \qquad A_3: \quad \sqrt{1/3}t_3,$

$A_4: \quad \text{double exponential}, \qquad A_5: \quad [0.5N(-3,1) + 0.5N(3,1)]/\sqrt{10}.$

The experiments are carried out using Fortran 77 and the optimization algorithm from the Fortran subroutine DBCOAH in the IMSL library. From Table 3, we can see that the size of $K_n$-test is somewhat conservative when $n = 200$. However, when the sample size $n$ increases to 400, the size is very close to the nominal significance level, for all selected levels. The power of this test depends on the shape of the alternative distribution. It increases as one moves from $A_1$ to $A_3$. This is reasonable since the difference between the shapes of $N(0, 1)$ and $A_3$ is more than that of the shapes of $N(0, 1)$ and $A_1$, and similarly for $A_2$ and $A_3$. At the alternative $A_4$, the power ranges between its values at the alternatives $A_2$ and $A_3$. The test reaches its highest power when the alternative is mixed-normal, that is, $A_5$. Some simulation results not reported here show that the empirical size and power of $K_n$ change a little when the true parameters $(a, \alpha_0, \alpha, \beta)$ are changed in the stationary region. These simulations indicate that the proposed test has satisfactory size and power behavior in finite samples, and should be useful in practice.

We next use model (3.6) to investigate the Hang Seng Index (HSI) in the Hong Kong stock market. Each period of two years from 1/6/1988–31/5/1996 is considered. The results are summarized in Table 4. In this table, the values in parentheses are the corresponding asymptotic standard deviations of the estimated parameters, LF is the value of log-likelihood function and QM is the portmanteau test statistic as in [20]. Both statistics, QM(6) and QM(12), suggest that this model fits the data adequately. All these estimators of the parameters $(\alpha, \beta)$ satisfy the finite fourth moment condition: $3\alpha^2 + 2\alpha\beta + \beta^2 < 1$. We use the statistic $K_n$ to test the hypothesis

TABLE 3
*Size and power of test statisitc $K_n$ for $\eta_i$ in* AR(1)–GARCH(1, 1) *models:* 1000 *replications*

| $\alpha$ $\eta_i \sim$ | $n = 200$ | | | | $n = 400$ | | |
|---|---|---|---|---|---|---|---|
| | **0.1** | **0.05** | **0.01** | | **0.1** | **0.05** | **0.01** |
| | | | | Size | | | |
| $N(0,1)$ | 0.089 | 0.041 | 0.006 | | 0.102 | 0.053 | 0.008 |
| | | | | Power | | | |
| $A_1$ | 0.171 | 0.086 | 0.021 | | 0.348 | 0.226 | 0.058 |
| $A_2$ | 0.309 | 0.180 | 0.056 | | 0.590 | 0.453 | 0.200 |
| $A_3$ | 0.570 | 0.434 | 0.201 | | 0.909 | 0.882 | 0.581 |
| $A_4$ | 0.407 | 0.247 | 0.060 | | 0.793 | 0.640 | 0.283 |
| $A_5$ | 1.00 | 1.00 | 1.00 | | 1.00 | 1.00 | 1.00 |

TABLE 4
*Empirical results for Hong Seng Index fitted* AR(1)–GARCH(1, 1) *models*

| **Periods** | $n$ | $a$ | $\alpha_0$ | $\alpha$ | $\beta$ | **LF** | **QM(6)** | **QM(12)** | $K_n$ |
|---|---|---|---|---|---|---|---|---|---|
| 1/6/88–31/5/90 | 493 | 0.242 | 0.083 | 0.223 | 0.772 | $-117.9$ | 5.6 | 12.3 | 6.44 |
| | | (0.055) | (0.027) | (0.048) | (0.031) | | | | |
| 1/6/90–31/5/92 | 495 | 0.186 | 0.526 | 0.203 | 0.442 | $-131.7$ | 2.9 | 7.6 | 4.67 |
| | | (0.056) | (0.162) | (0.070) | (0.145) | | | | |
| 1/6/92–31/5/94 | 498 | 0.117 | 0.322 | 0.242 | 0.664 | 25.5 | 6.3 | 14.5 | 1.74 |
| | | (0.050) | (0.108) | (0.057) | (0.068) | | | | |
| 1/6/94–31/5/96 | 497 | 0.128 | 0.053 | 0.069 | 0.900 | $-95.6$ | 7.9 | 12.9 | 2.16 |
| | | (0.048) | (0.029) | (0.025) | (0.034) | | | | |

$H_0: \eta_i \sim N(0,1)$. Upon comparing the values of $K_n$ in Table 4 with the critical values obtained in Table 4, one observes that the test $K_n$ rejects $H_0$ at significance level 0.01 during the first two periods, while it accepts $H_0$ at significance level 0.10 during other two periods. Note that the Tian An Men Square event occurred during the summer of 1989. This event had a big impact on the Hong Kong stock market such that HSI has some nonnormal factors during this period. The event affected the Hong Kong stock market until 1992 and after that, its effect gradually disappeared. No other exciting event occurred and the HSI was rather quiet during 1992–1996. Our findings here seem to match this real circumstance very well.

**4. Some general results.** This section contains some general results useful for obtaining various approximations in hereroscedastic nonlinear time series. The first result is an extension of Theorem 1.1 in [17]. To state these extensions, let $(\eta_{ni}, \gamma_{ni}, \tau_{ni}, \xi_{ni}), 1 \leq i \leq n$, be an array of 4-tuple r.v.'s defined on a probability space such that $\{\eta_{ni}, 1 \leq i \leq n\}$ are i.i.d. according to a d.f. $H$, and $\eta_{ni}$ is independent of $(\gamma_{ni}, \tau_{ni}, \xi_{ni}), 1 \leq i \leq n$. Furthermore, let $\{\mathcal{A}_{ni}\}$ be an array of sub $\sigma$-fields such that $\mathcal{A}_{ni} \subset \mathcal{A}_{ni+1}, 1 \leq i \leq n, n \geq 1$; $(\gamma_{n1}, \delta_{n1}, \xi_{n1})$ is $\mathcal{A}_{n1}$ measurable, the r.v.'s $\{\eta_{n1}, \ldots, \eta_{nj-1}; (\gamma_{ni}, \tau_{ni}, \xi_{ni}), 1 \leq i \leq j\}$ are $\mathcal{A}_{nj}$-measurable, $2 \leq j \leq n$; and $\eta_{nj}$ is in dependent of $\mathcal{A}_{nj}, 1 \leq j \leq n$. Define, for $x \in \mathbb{R}$,

$$\tilde{V}_n(x) := n^{-1/2} \sum_{i=1}^n \gamma_{ni} I(\eta_{ni} \leq x + x\tau_{ni} + \xi_{ni}),$$

$$\tilde{J}_n(x) := n^{-1/2} \sum_{i=1}^n \gamma_{ni} H(x + x\tau_{ni} + \xi_{ni}),$$

$$V_n(x) := n^{-1/2} \sum_{i=1}^n \gamma_{ni} I(\eta_{ni} \leq x + \xi_{ni}), \qquad J_n(x) := n^{-1/2} \sum_{i=1}^n \gamma_{ni} H(x + \xi_{ni}),$$

$$V_n^*(x) := n^{-1/2} \sum_{i=1}^n \gamma_{ni} I(\eta_{ni} \leq x), \qquad J_n^*(x) := n^{-1/2} \sum_{i=1}^n \gamma_{ni} H(x),$$

$$\tilde{U}_n(x) := \tilde{V}_n(x) - \tilde{J}_n(x), \qquad U_n(x) := V_n(x) - J_n(x),$$

$$U_n^*(x) := V_n^*(x) - J_n^*(x).$$

Next, we introduce some assumptions:

(4.1) $\quad H$ has a.e. positive density $h$ with $\|h\|_\infty := \sup_{x \in \mathbb{R}} h(x) < \infty$.

(4.2) $$\int |x| h(x)\, dx < \infty.$$

(4.3) $\quad$ (a) $\quad n^{-1} \sum_{i=1}^n \gamma_{ni}^2 = O_p(1)$, $\quad$ (b) $\quad \max_{1 \leq i \leq n} n^{-1/2} |\gamma_{ni}| = o_p(1)$.

(4.4) $\quad$ (a) $\quad \max_{1 \leq i \leq n} |\xi_{ni}| = o_p(1)$, $\quad$ (b) $\quad \max_{1 \leq i \leq n} |\tau_{ni}| = o_p(1)$.

A close inspection of the proof of Theorem 1.1 of [17] and the discussion on the bottom of page 544 there shows that if (4.1), (4.3), and (4.4)(a) hold, then

(4.5) $$\sup_{x \in \mathbb{R}} |U_n(x) - U_n^*(x)| = o_p(1).$$

This result has played a pivotal role in the development of a unified approach to the asymptotic distribution analysis of numerous inference procedures in homoscedastic time series models under minimal moment assumptions; see [15, 16]. We now state an analog of (4.5) suitable for heteroscedastic models.

THEOREM 4.1. *Under the above setup and under assumptions* (4.1)–(4.4),

(4.6) $$\sup_{x\in\mathbb{R}} |\tilde{U}_n(x) - U_n^*(x)| = o_p(1).$$

The proof of (4.5) used a chaining argument with respect to the pseudo-metric $d_b(x,y) := \sup_{|z|\le b} |H(x+z) - H(y+z)|^{1/2}, x,y \in \mathbb{R}, b > 0$, appropriate for the location problem. If we let $\mathcal{N}(\delta,b)$ denote the cardinality of the minimal $\delta$-net of $(\mathbb{R}, d_b)$, then the crucial condition needed for this chaining argument is that, for some $0 < b_0 < 1$, $\int_0^1 \{\ln \mathcal{N}(u,b_0)\}^{1/2} du < \infty$. On page 544 of [17] it is argued that (4.1) implies this condition.

The analog of the above metric that works in the current setup, that is, for the location-scale setup, is defined as follows. Let $|z| := |z_1| \vee |z_2|$, for any $z := (z_1, z_2)' \in \mathbb{R}^2$, and define

$$\rho_b(x,y) = \sup_{|z|\le b} |H(x(1+z_1) + z_2) - H(y(1+z_1) + z_2)|^{1/2}, \qquad x,y \in \mathbb{R}, b > 0.$$

Let $N(\delta, b)$ be the cardinality of the minimal $\delta$-net of $(\mathbb{R}, \rho_b)$ and let

$$I(b) := \int_0^1 \{\ln N(u,b)\}^{1/2} du.$$

We shall now show that (4.1) and (4.2) imply

(4.7) $$I(b) < \infty \qquad \forall\, 0 \le b < 1.$$

Fix a $0 < \delta, b < 1$. Since $H$ is continuous, there exists $0 < M_{1\delta}, M_{2\delta} < \infty$ such that

$$H(-M_{1\delta}) = \delta^2/2 \quad \text{and} \quad 1 - H(M_{2\delta}) = \delta^2/2.$$

When $x, y \in (-\infty, -(b+M_{1\delta})/(1-b)]$,

$$\rho_b(x,y) \le |H(x(1-b) + b) + H(y(1-b) + b)|^{1/2} \le \{2H(-M_{1\delta})\}^{1/2} < \delta.$$

For $x, y \in [(b+M_{2\delta})/(1-b), \infty)$, using the monotonicity of $H$,

$$\rho_b(x,y) = \sup_{|z|\le b} |H(x(1+z_1) + z_2) - H(y(1+z_1) + z_2)|^{1/2}$$
$$\le [|1 - H(x(1-b) - b)| + |1 - H(y(1-b) - b)|]^{1/2}$$
$$\le \{2(1 - H(M_{2\delta}))\}^{1/2} < \delta.$$

Partition the interval $[-(b+M_{1\delta})/(1-b), (b+M_{2\delta})/(1-b)]$ as $-(b+M_{1\delta})/(1-b) = x_1 < x_2 < \cdots < x_{N_\delta} = (b+M_{2\delta})/(1-b)$ with $x_k - x_{k-1} = \delta^2/(2\|h\|_\infty)$ and $N_\delta = 2(2b + M_{1\delta} + M_{2\delta})\|h\|_\infty/[(1-b)\delta^2]$, where $[z]$ is the integer part of $z \in \mathbb{R}$. For each $x, y \in (x_{k-1}, x_k]$, using (4.1), we obtain

$$\rho_b(x,y) = \sup_{|z| \leq b} |H(x(1+z_1) + z_2) - H(y(1+z_1) + z_2)|^{1/2}$$

$$\leq [\|h\|_\infty |x-y|(1+b)]^{1/2} < \delta \qquad (\text{because } b < 1).$$

With $\mu := E|\eta_{n1}|$, by (4.2) and the Markov inequality, $M_{1\delta} + M_{2\delta} \leq 4\mu/\delta^2$, because

$$\frac{\delta^2}{2} = H(-M_{1\delta}) \leq P(|\eta_{n1}| \geq M_{1\delta}) \leq \frac{\mu}{M_{1\delta}},$$

$$\frac{\delta^2}{2} = 1 - H(M_{2\delta}) \leq P(|\eta_{n1}| \geq M_{2\delta}) \leq \frac{\mu}{M_{2\delta}}.$$

Hence, using $\delta < 1$, we have

$$N(\delta, b) \leq 2 + N_\delta \leq 2 + \frac{2(2b + M_{1\delta} + M_{2\delta})\|h\|_\infty}{(1-b)\delta^2}$$

$$\leq 2 + \frac{2(2b\delta^2 + 4\mu)\|h\|_\infty}{(1-b)\delta^4} \leq \frac{2(1-b)\delta^4 + 2(2b\delta^2 + 4\mu)\|h\|_\infty}{(1-b)\delta^4}$$

$$\leq \frac{2(1-b) + 2(2b + 4\mu)\|h\|_\infty}{(1-b)\delta^4} = \frac{C_b}{\delta^4}.$$

Note that $C_b := 2 + [(4b + 8\mu)\|h\|_\infty/(1-b)]$ is an increasing function of $b$, with $\sup_{0 \leq b \leq b_1} C_b = C_{b_1} < \infty$, for any $b_1 < 1$. This proves (4.7).

Thus, one can now repeat the arguments in [17] verbatim using the metric $\rho_b$ instead of $d_b$. No details are given here. However, we mention one typo that appears in the statement of Proposition 2.1 of [17] on page 552. The event $[\max_i |\gamma_{ni}| \leq \delta/(1 + \ln N(\delta, b))^{1/2}]$ at (2.5) there has $n^{1/2}$ missing. It should be $[\max_i |\gamma_{ni}| \leq n^{1/2}\delta/(1 + \ln N(\delta, b))^{1/2}]$.

The above Theorem 4.1 is not enough to cover the cases where the weights $\gamma_{ni}$ and the disturbances $\tau_{ni}, \xi_{ni}$ are functions of certain parameters and where one desires to obtain various approximations uniformly in these parameters, as is needed in the previous sections. The next result gives the needed extension of this theorem to cover these cases. Accordingly, let $m \geq 1$ be a fixed integer, $l_{ni}, v_{ni}, u_{ni}$ be measurable functions from $\mathbb{R}^m$ to $\mathbb{R}$ such that, for every $t \in \mathbb{R}^m$, $(l_{ni}(t), v_{ni}(t), u_{ni}(t))$ are independent of $\eta_{ni}$, and $\mathcal{A}_{ni}$-measurable, for each $1 \leq i \leq n$. Let, for $x \in \mathbb{R}, t \in \mathbb{R}^m$,

$$\mathcal{V}(x,t) := n^{-1/2} \sum_{i=1}^n l_{ni}(t) I(\eta_{ni} \leq x + xv_{ni}(t) + u_{ni}(t)),$$

$$\mathcal{J}(x,t) := n^{-1/2} \sum_{i=1}^n l_{ni}(t) H(x + xv_{ni}(t) + u_{ni}(t)),$$

(4.8)
$$\tilde{\mathcal{U}}(x,t) := \mathcal{V}(x,t) - \mathcal{J}(x,t),$$
$$\mathcal{U}^*(x,t) := n^{-1/2} \sum_{i=1}^{n} l_{ni}(t)[I(\eta_{ni} \leq x) - H(x)].$$

To state the needed result, we need the following assumptions. For each $t \in \mathbb{R}^m$,

(4.9) $$n^{-1} \sum_{i=1}^{n} l_{ni}^2(t) = O_p(1), \qquad \max_{1 \leq i \leq n} n^{-1/2} |l_{ni}(t)| = o_p(1),$$

(4.10) $$\max_{1 \leq i \leq n} \{|v_{ni}(t)| + |u_{ni}(t)|\} = o_p(1),$$

(4.11) $$n^{-1/2} \sum_{i=1}^{n} |l_{ni}(t)|[|v_{ni}(t)| + |u_{ni}(t)|] = O_p(1).$$

$\forall \epsilon > 0, \exists \delta > 0$, and an $n_1 \ni \forall 0 < b < \infty, \forall \|s\| \leq b, \forall n > n_1$,

(4.12)
$$P\left(n^{-1/2} \sum_{i=1}^{n} |l_{ni}(s)| \left\{ \sup_{\|t-s\|<\delta} |v_{ni}(t) - v_{ni}(s)| \right.\right.$$
$$\left.\left. + \sup_{\|t-s\|<\delta} |u_{ni}(t) - u_{ni}(s)| \right\} \leq \epsilon \right) > 1 - \epsilon,$$

(4.13) $$P\left(\sup_{\|t-s\|<\delta} n^{-1/2} \sum_{i=1}^{n} |l_{ni}(t) - l_{ni}(s)| \leq \epsilon \right) > 1 - \epsilon.$$

The following lemma gives the needed result.

LEMMA 4.1. *Under the above setup and under the assumptions* (4.1) *and* (4.9)–(4.13), *for every* $0 < b < \infty$, $\sup_{x \in \mathbb{R}, \|t\| \leq b} |\tilde{\mathcal{U}}(x,t) - \mathcal{U}^*(x,t)| = o_p(1)$.

PROOF. The proof of this lemma is similar to that of Lemma 8.3.2 in [16], with the proviso that one uses the above Theorem 4.1 instead of Theorem 2.2.5 of [16] whenever needed. □

PROOF OF THEOREM 2.1. Proof of Theorem 2.1 follows from Lemma 4.1 applied as follows. Fix a $1 \leq j \leq r$, and let $W_{22j,i}$ denote the $j$th component of the vector $W_{22,i}$. Now in the above setup take

$$u_{ni}(t) := \frac{\mu_i(\theta + n^{-1/2}t) - \mu_i(\theta)}{\sqrt{h_i(\theta)}},$$

$$(4.14) \quad v_{ni}(t) := \frac{\sqrt{h_i(\theta + n^{-1/2}t)} - \sqrt{h_i(\theta)}}{\sqrt{h_i(\theta)}},$$

$$\ell_{ni}(t) := W_{22j,i}(\theta + n^{-1/2}t), \qquad 1 \le i \le n, t \in \mathbb{R}^{q+r}, m = q+r.$$

Rewrite in terms of the above notation:

$$I(\eta_i(\theta + n^{-1/2}t) \le x) \equiv I(\eta_i \le x + xv_{ni}(t) + u_{ni}(t)),$$

$$K_n(x, \theta + n^{-1/2}t) = \tilde{\mathcal{U}}(x,t)$$
$$+ n^{-1/2} \sum_{i=1}^{n} \ell_{ni}(t)[F_0(x + xv_{ni}(t) + u_{ni}(t)) - F_0(x)],$$

where in $\tilde{\mathcal{U}}$, the d.f. $H$ is replaced by $F_0$.

The rest of the proof consists of verifying the conditions (4.9)–(4.13) for the entities given at (4.14), the details of which are similar to those appearing in the proof of Lemma 8.3.1 in [16], hence omitted for the sake of brevity. □

PROOF OF THEOREM 2.2. The result follows by use of Le Cam's third lemma (cf. [30], Theorem 3.10.7), as soon as we verify the contiguity of $P_{1n}$ to $P_{0n}$. The sequence of probability measures $\{P_{1n}\}$ here depends on $n$ through the change in direction of the error distribution, while in [8], analogous alternative sequences depend on $n$ through the change in a Euclidean parameter. Thus, their results of LAN are inapplicable.

Now, let $f_c = f + c(\tilde{f} - f)$ and $\zeta = (\tilde{f} - f)/\sqrt{f}$, $|c| < 1$ and $\tilde{f}$ be the density of $\tilde{F}$. Assume that $0 < \sigma^2 \equiv \int \zeta^2(x) f(x) \, dx < \infty$. Let $\lambda_n(c) := \sum_{i=1}^{n} \log(f_c(\eta_i)/f(\eta_i))$ denote the log-LR of $P_{1n}^c$ to $P_{0n}$. Verify that $c^{-2} \int [\sqrt{f_c(x)} - \sqrt{f(x)} - \frac{1}{2}c\zeta(x)\sqrt{f(x)}]^2 \, dx \to 0$, as $c \to 0$. From the proof of Theorem 2.1 in [26], under $P_{0n}$ it follows that

$$\lambda_n(n^{-1/2}\delta) = \delta n^{-1/2} \sum_{i=1}^{n} \zeta_i(\eta_i(\theta)) - \tfrac{1}{2}\delta^2 \sigma^2 + o_p(1).$$

By the central limit theorem, it now readily follows that, under $P_{0n}$, $n^{-1/2} \times \sum_{i=1}^{n} \zeta_i(\eta_i(\theta)) \xrightarrow{d} N(0, \sigma^2)$, which in turn implies that $P_{1n}$ is contiguous to $P_{0n}$. □

**Acknowledgments.** The authors thank a referee, an Associate Editor and the Editor J. Fan for their helpful comments.

## REFERENCES

[1] BERKES, I. and HORVÁTH, L. (2003). Limit results for the empirical process of squared residuals in GARCH models. *Stochastic Process. Appl.* **105** 271–298. MR1978657

[2] BOLDIN, M. V. (1982). Estimation of the distribution of noise in an autoregressive scheme. *Theory Probab. Appl.* **27** 866–871. MR0681475

[3] BOLDIN, M. V. (1989). On hypothesis testing in a moving average scheme by the Kolmogorov–Smirnov and the $\omega^2$ tests. *Theory Probab. Appl.* **34** 699–704. MR1036715

[4] BOLDIN, M. V. (1998). On residual empirical distribution functions in ARCH models with applications to testing and estimation. *Mitt. Math. Sem. Giessen* No. 235 49–66. MR1661093

[5] BOLDIN, M. V. (2002). On sequential residual empirical processes in heteroscedastic time series. *Math. Methods Statist.* **11** 453–464. MR1979744

[6] BOUGEROL, P. and PICARD, N. M. (1992). Stationarity of GARCH processes and of some nonnegative time series. *J. Econometrics* **52** 115–127. MR1165646

[7] D'AGOSTINO, R. B. and STEPHENS, M. A. (1986). *Goodness-of-Fit Techniques*. Dekker, New York. MR0874534

[8] DROST, F. C., KLAASSEN, C. A. J. and WERKER, B. J. M. (1997). Adaptive estimation in time-series models. *Ann. Statist.* **25** 786–817. MR1439324

[9] DURBIN, J. (1973). *Distribution Theory for Tests Based on the Sample Distribution Function*. SIAM, Philadelphia. MR0305507

[10] FRANCQ, C. and ZAKOÏAN, J.-M. (2004). Maximum likelihood estimation of pure GARCH and ARMA–GARCH processes. *Bernoulli* **10** 605–637. MR2076065

[11] HORVÁTH, L., KOKOSZKA, P. and TEYSSIÈRE, G. (2001). Empirical process of the squared residuals of an ARCH sequence. *Ann. Statist.* **29** 445–469. MR1863965

[12] HORVÁTH, L., KOKOSZKA, P. and TEYSSIÈRE, G. (2004). Bootstrap misspecification tests for ARCH based on the empirical process of squared residuals. *J. Stat. Comput. Simul.* **74** 469–485. MR2073226

[13] KOUL, H. L. (1991). A weak convergence result useful in robust autoregression. *J. Statist. Plann. Inference* **29** 291–308. MR1144174

[14] KOUL, H. L. (1992). *Weighted Empiricals and Linear Models*. IMS, Hayward, CA. MR1218395

[15] KOUL, H. L. (1996). Asymptotics of some estimators and sequential residual empiricals in nonlinear time series. *Ann. Statist.* **24** 380–404. MR1389897

[16] KOUL, H. L. (2002). *Weighted Empirical Processes in Dynamic Nonlinear Models*, 2nd ed. *Lecture Notes in Statist.* **166**. Springer, New York. MR1911855

[17] KOUL, H. L. and OSSIANDER, M. (1994). Weak convergence of randomly weighted dependent residual empiricals with applications to autoregression. *Ann. Statist.* **22** 540–562. MR1272098

[18] KULPERGER, R. and YU, H. (2005). High moment partial sum processes of residuals in GARCH models and their applications. *Ann. Statist.* **33** 2395–2422. MR2211090

[19] LE CAM, L. (1986). *Asymptotic Methods in Statistical Decision Theory*. Springer, New York. MR0856411

[20] LI, W. K. and MAK, T. K. (1994). On the squared residual autocorrelations in non-linear time series with conditional heteroskedasticity. *J. Time Ser. Anal.* **15** 627–636. MR1312326

[21] LING, S. (1998). Weak convergence of the sequential empirical processes of residuals in nonstationary autoregressive models. *Ann. Statist.* **26** 741–754. MR1626028

[22] LING, S. (1999). On the probabilistic properties of a double threshold ARMA conditional heteroskedastic model. *J. Appl. Probab.* **36** 688–705. MR1737046

[23] Ling, S. (2003). Adaptive estimators and tests of stationary and non-stationary short- and long- memory ARFIMA–GARCH models. *J. Amer. Statist. Assoc.* **98** 955–967. [MR2041484](MR2041484)

[24] Ling, S. (2005). Self-weighted and local quasi-maximum likelihood estimators for ARMA–GARCH /IGARCH models. *J. Eonometrics*. To appear.

[25] Ling, S. and Li, W. K. (1997). On fractionally integrated autoregressive moving-average time series models with conditional heteroskedasticity. *J. Amer. Statist. Assoc.* **92** 1184–1194. [MR1482150](MR1482150)

[26] Ling, S. and McAleer, M. (2003). On adaptive estimation in nonstationary ARMA models with GARCH errors. *Ann. Statist.* **31** 642–674. [MR1983545](MR1983545)

[27] Ling, S. and McAleer, M. (2004). Regression quantiles for unstable autoregression models. *J. Multivariate Anal.* **89** 304–328. [MR2063636](MR2063636)

[28] Loynes, R. M. (1980). The empirical distribution function of residuals from generalized regression. *Ann. Statist.* **8** 285–298. [MR0560730](MR0560730)

[29] Tsay, R. S. (2005). *Analysis of Financial Time Series*, 2nd ed. Wiley, New York. [MR2162112](MR2162112)

[30] van der Vaart, A. W. and Wellner, J. A. (1996). *Weak Convergence and Empirical Processes*: *With Applications to Statistics*. Springer, New York. [MR1385671](MR1385671)

Department of Statistics and Probability  
Michigan State University  
East Lansing, Michigan 48824-1027  
USA  
E-mail: [koul@stt.msu.edu](koul@stt.msu.edu)

Department of Mathematics  
Hong Kong University  
of Science and Technology  
Clear Water Bay  
Hong Kong  
E-mail: [maling@ust.hk](maling@ust.hk)